\documentclass[12pt]{amsart}
\usepackage{amssymb,amscd}
\usepackage[mathscr]{eucal}

\input epsf
\newtheorem{thm}{Theorem}[section]
\newtheorem{lem}[thm]{Lemma}

\numberwithin{equation}{subsection}
\begin{document}

\title{Ergodicity of Mapping Class Group Actions 
on Representation Varieties, II. Surfaces with Boundary}
\author{Doug Pickrell and Eugene Z. Xia}
\address{
Department of Mathematics\\ 
University of Arizona\\
Tucson, AZ 85721 ({\it Pickrell})
}
\address{
National Center for Theoretical Sciences\\
Third General Building \\
National Tsing Hua University\\
No. 101, Sec 2, Kuang Fu Road, \\
Hsinchu, Taiwan 30043, Taiwan R.O.C.
({\it Xia}) \bigskip
}
\email{
pickrell@math.arizona.edu {\it (Pickrell)}, 
xia@math.umass.edu {\it (Xia)}}
\date{\today}

\renewcommand{\o}{\operatorname}
\newcommand{\C}{{\mathbb C}}
\newcommand{\CC}{{\mathcal C}}
\renewcommand{\P}{{\mathbb P}}
\newcommand{\bH}{{\mathbb H}}
\renewcommand{\H}{\o{H}}
\newcommand{\R}{{\mathbb R}}
\newcommand{\Z}{{\mathbb Z}}
\newcommand{\Id}{{\mathbb I}}
\newcommand{\Os}{{\mathcal O}}
\newcommand{\zz}{{\mathcal Z}}
\newcommand{\z}{{\zz}^1(M)}
\newcommand{\M}{{\mathcal M}}
\newcommand{\la}{\langle}
\newcommand{\ra}{\rangle}
\newcommand{\h}{{\mathcal H}^3_{\R}}
\newcommand{\hH}{{\mathscr H}}
\newcommand{\hi}{h^{-1}}
\renewcommand{\Im}{\o{Im}}
\renewcommand{\Re}{\o{Re}}
\newcommand{\Hol}{\o{Hol}(E)}
\newcommand{\Holu}{\o{Hol}_u(E)}
\newcommand{\hol}{\o{hol}}
\newcommand{\her}{\o{Her}}
\newcommand{\Her}{\her(E)}
\newcommand{\fc}{{\mathcal F}(E)} 
\newcommand{\fcu}{{\mathcal F}_u(E)} 
\newcommand{\fcuu}{{\mathcal F}_u(E)} 
\newcommand{\fcl}{{\mathcal F}_l(E)} 
\newcommand{\g}{{\mathfrak g}}
\renewcommand{\k}{{\mathfrak k}}
\newcommand{\gl}{{\mathfrak{gl}}}
\newcommand{\su}{{\mathfrak{su}}}
\newcommand{\ff}{{\mathfrak F}}	
\newcommand{\ga}{{\mathcal G}(E)} 
\newcommand{\gau}{{\mathcal G}_u(E)} 
\newcommand{\gal}{{\mathcal G}_l(E)} 
\newcommand{\Hig}{\o{Higgs}(E)}
\newcommand{\Hom}{\o{Hom}}
\newcommand{\hpg}{\Hom(\pi,G)/G}
\newcommand{\hppg}{\Hom(\pi,G)\slashslash G}
\newcommand{\hmg}{\Hom(\pi,G)}
\newcommand{\Un}{\o{U}(n)}
\newcommand{\Uo}{\o{U}(1)}
\newcommand{\GLn}{\o{GL}(n,\C)}
\newcommand{\GL}{\o{GL}}
\newcommand{\PUoo}{\o{PU}(1,1)}
\newcommand{\PSLtR}{\o{PSL}(2,\R)}
\newcommand{\PSLnR}{\o{PSL}(n,\R)}
\newcommand{\SLnC}{\o{SL}(n,\C)}
\newcommand{\PSLnC}{\o{PSL}(n,\C)}
\newcommand{\SOto}{\o{SO}(2,1)}
\newcommand{\rk}{\o{rank}}
\newcommand{\Map}{\o{Map}}
\newcommand{\db}{\overline{\partial}}
\renewcommand{\d}{\partial}
\newcommand\ka{K\"ahler~}
\newcommand\tM{\tilde{M}}
\newcommand\tS{\tilde{S}}
\newcommand\td{\tilde{d}}
\newcommand\tg{\tilde{g}}
\newcommand\tT{\tilde{T}}
\renewcommand\th{\tilde{h}}
\newcommand\ca{{\mathcal A}}
\newcommand\cb{{\mathcal B}}
\newcommand\cc{{\mathcal C}}
\newcommand\Iso{{\o{Iso}}}
\newcommand\Aut{{\o{Aut}}}
\newcommand\Ad{{\o{Ad}}}
\newcommand\Obj{{\o{Obj}}}
\newcommand\hk{{hyperk\"ahler~}}
\newcommand\slashslash{{/\hspace{-3pt}/}}
\newcommand\Er{E_{\rho}}
\begin{abstract}
The mapping class
group of a compact oriented surface of genus greater than
one with boundary acts
ergodically on connected components of the
representation variety corresponding to a connected
compact Lie group, for every choice of conjugacy class
boundary condition.
\end{abstract}
\maketitle

\section{Introduction and Results}
Let $\Sigma$ be a compact oriented surface 
with $n+1$ boundary components (circles) $(c_0,...,c_n)$ and $K$ a compact
connected Lie group of dimension $d$ and rank $r$.
Denote by $\Gamma$ the mapping class group of $\Sigma$
that fixes the $c_j's$.  
Fix a set of conjugacy classes $\{C_j\}_{j=0}^n$, one for each
boundary component and let 
$$
C = \prod_{j=0}^n C_j.
$$
Then the representation variety is
$$
\Hom_C(\pi_1(\Sigma), K) = \{\rho : \rho(c_j) \in C_j \}.
$$
The group $K$ acts on $\Hom_C(\pi_1(\Sigma),K)$ by conjugation and
the moduli space is the quotient 
$$
\Hom_C(\pi_1(\Sigma), K)/conj(K).
$$

\begin{thm} ~\label{main}
Fix $p>1$.
With respect to the Lebesgue class, $\Gamma$ acts ergodically on the
connected components of $\Hom_C(\pi_1(\Sigma), K)/conj(K)$.
\end{thm}
When $n=0$, this theorem was proved in \cite{PX1}.

\section{The representation variety and the mapping class group}
With slight modifications,
the basic notations are as in \cite{PX1}.
Let $K' = [K,K]$ be the semi-simple part of $K$ and $dg', dg$
their respective Haar measures.  Denote by $\k, \k'$ 
the Lie algebras of $K,K'$, respectively.

As $\Sigma$ has boundary, by adding a 
possible extra boundary component and a fixed central element
for that boundary component, one may assume $K'$ is simply
connected without losing generality.

\bigskip
\bigskip

\centerline{\epsfysize=4in            \epsffile{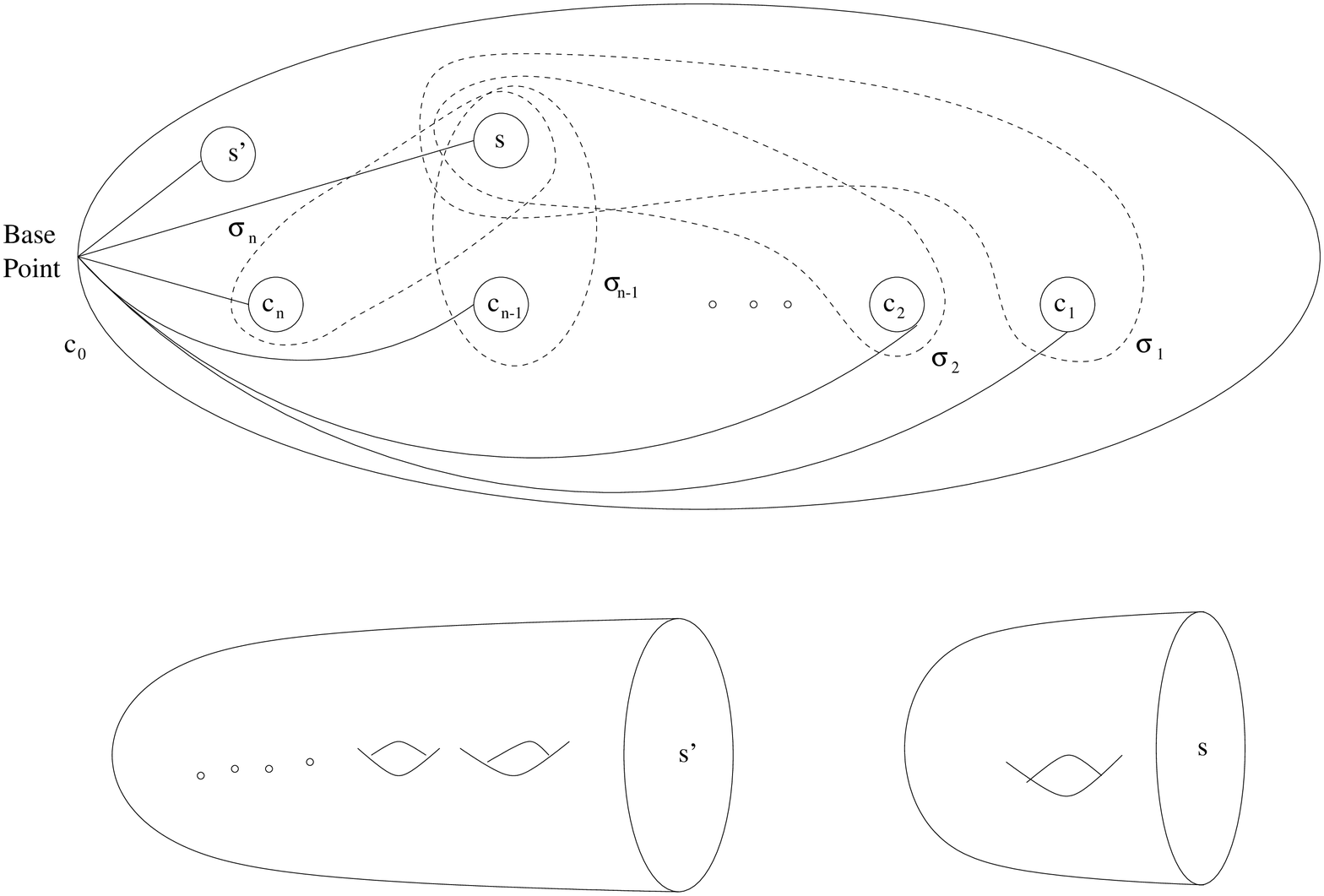}}
\centerline{F{\sc igure} 1: The surface $\Sigma$.}

\

Fix the base point on $c_0$.
Then $\pi_1(\Sigma)$ has the following presentation:
$$
\la x_i,y_i, c_j \ \ |\ \prod_{i=1}^p [x_i,y_i] \prod_{j=0}^n c_j \ra.
$$
In terms of Figure 1, $x_1,y_1$ belong to the one-holed torus
bounded by $s$ while for $i > 1$, $x_i,y_i$ belong to the
one-holed genus $(p-1)$ surface bounded by $s'$.
All the boundary circles in Figure 1 are oriented 
counter-clockwise with the exception of $c_0$.
The following notions are fixed through out this paper.
By a slight abuse of notation,
the $c_j's$ will also denote elements in $K$.
Let $g = (g_1,...,g_p), h = (h_1,...,h_p) \in K^{p}$
and $c = (c_0,...,c_n) \in K^{n+1}$, although in many cases
$c$ will denote elements in $C$.
Define maps
$$
R_{l,k}(g,h) = \prod_{i=l}^k [g_i,h_i],
$$
$$
S_{l,k}(c) = \prod_{i=l}^k c_j.
$$
By Remark (2.1.5) in \cite{PX1}, the image of the map
$R_{l,k} : K^{2(k-l)} \longrightarrow K'$ has full
measure in $K'$.  Since the image of $R_{l,k}$ is compact,
$R_{l,k}$ is actually onto $K'$.

The mapping class group $\Gamma$ is the group of connected components
of the diffeomorphism group of $\Sigma$ fixing the boundary
$\partial \Sigma$.   The group $\Gamma$ is generated by the 
Dehn twists.  In particular, $\Gamma$ contains the
twists along $x_i,y_i$ and $\sigma_j$ portrayed
in Figure 1.  Let $\Gamma_0$ be the subgroup generated
by the Dehn twists along the $x_i$'s and $y_i$'s
and $\Gamma_j$ be the subgroup generated by $\Gamma_0 \cup \{\sigma_j\}$.

\section{Ergodicity}
Let $p > 1$.
Fix conjugacy classes $C_j$ for each, $0 \le j \le n$.
Then the representation variety $\Hom_C(\Sigma,K)$ is
$$
\{ (g,h,c) \in  K^{2p} \times C: 
R_{1,p}(g,h) S_{0,n}(c) = I\}.
$$
Assume $\Hom_C(\Sigma,K)$ is non-empty which means that
$S_{0,n}(c) \in K'$.
Each conjugacy class $C_j$ possess a $K$-invariant canonical 
measure.  From these measures on the conjugacy classes and
the measure $dg$ on $K$,
the variety $\Hom_C(\Sigma,K)$ inherits a Lebesgue measure whose
push-forward measure on $\Hom_C(\pi_1(\Sigma),K)/conj(K)$ is in
the same measure class as $\mu$.

Fix $j$ and for each $l \neq j$, fix an element $c_l \in C_l$.
Let $k = S_{l+1,n}(c)$, $s = R_{1,1}(g,h)$
and $s' = R_{2,p}(g,h)$.
Then the action of $\sigma_j$ fixes $c_l$ for all $j \neq l$ (See Figure 1)
and $\sigma_j (c_j) = (s k c_j k^{-1}) c_j (s k c_j k^{-1})^{-1}$.
Let 
$$
M = \{(g,h,c_j): R_{1,p}(g, h) S_{0,n}(c)= I \} \subset 
\Hom_C(\pi_1(\Sigma),K).
$$
For each fixed $c_j \in C_j$, let 
$$
N(c_j) = \{(g,h): R_{1,p}(g, h) S_{0,n}(c)= I \} \subset M.
$$  
Then both
$M$ and $N(c_j)$ inherit measures from  $\Hom_C(\Sigma,K)$.  Indeed,
$M$ and $N(c_j)$ are themselves both representation varieties on a surface
with two boundary components.  In the case of $M$, one boundary
is associated with a fixed group element and the other a conjugacy
class while both boundary components of $N(c_j)$ are associated
with group elements.  In particular, as subspaces of $\Hom_C(\Sigma,K)$, 
$M$ is $\Gamma_j$-invariant while $N(c_j)$ is $\Gamma_0$-invariant.

Let $f \in L^2(M)$ be $\Gamma_j$-invariant.  
The goal is to show that $f$ is constant a.e.
The following theorem is from \cite{PX1}:
\begin{thm}
The group $\Gamma_0$ acts ergodically on $N(c_j)$.
\end{thm}
This shows that $f$ is
constant on $N(c_j)$ for a.e. $c_j \in C_j$.
Hence $f$ is a characteristic function
on $C_j$.  It remains to show that $f$ is independent
of $c_j$.
Consider the function
$$
\phi : M \longrightarrow C_j \times C_j :
$$
$$
\ \ \ \ \phi(g,h,c_j) = (c_j,  (s k c_j k^{-1}) c_j (s k c_j k^{-1})^{-1}).
$$
Since $f$ is $\Gamma_j$-invariant, to show that $f$ is independent of $c_j$,
it is sufficient to show that the relation defined by $Im(\phi)$ is 
transitive; since the equality holds
in an a.e. sense, it is sufficient to show that

\begin{lem}~\label{lem:main}
Let $pr_1 : C_j \times C_j \longrightarrow C_j$ be the projection
to the first coordinate.  Then $pr_1(Interior(Im(\phi))) = C_j$.
\end{lem}
This Lemma implies that for each $c_j \in C_j$, there
are open sets $U_{c_j}, V_{c_j} \in C_j$ such that 
$c_j \in U_{c_j}$ and $U_{c_j} \times V_{c_j} \subset Im(\phi)$.
Since $\phi_1^* f = \phi_2^* f$, a.e., it follows that
$f$ is constant on $U_{c_j}$, a.e..  Since $pr_1(Interior(Im(\phi))) = C_j$
which is connected, $f$ is independent of $c_j$.
\begin{proof}
Since the $c_l$'s are all fixed for $l \neq j$, $k$ is fixed.
Since
$$
R_{1,1} :  K^2 \longrightarrow K',
$$
$$
R_{2,p} :  K^{2(p-1)} \longrightarrow K'
$$
are onto, 
one may choose regular values $s,s'$ respectively of these two
maps $R_{1,1}$ and $R_{2,p}$ such that $s s' S_{0,n}(c) = I$.
Let $(g,h,c_j) \in M$ such that $R_{1,1}(g,h) = s$.
Then $(g,h,c_j)$ is a smooth point 
(The representation corresponding to 
$(g,h,c_j)$ is irreducible) 
and the map $\phi_1$ is regular at 
$(g,h,c_j)$.
It will then be
sufficient to fix $c_j$ and show that the map 
$$
\phi_2 : N(c_j) \longrightarrow C_j :
$$
$$
\ \ \ \phi_2(g,h) =  (s k c_j k^{-1}) c_j (s k c_j k^{-1})^{-1}
$$
is regular at $(g,h)$.
Again, since both $s, s'$ are regular values
for the maps
$$
R_{1,1} :  K^2 \longrightarrow K',
$$
$$
R_{2,p} :  K^{2(p-1)} \longrightarrow K',
$$
$s$ is a regular value for the map
$$
R_{1,1} : N(c_j) \longrightarrow K'.
$$
This implies that the differential $d\phi_2|_{(g,h)}$
is surjective.
\end{proof}
Lemma~\ref{lem:main} shows that $f$ is independent of the variable
$c_j$ or the $\Gamma_j$-action is ergodic on $M$.  
If $F \in L^2(\Hom_C(\Sigma,K))$ is $\Gamma$-invariant,
then again by the main result in \cite{PX1}, $F$ is a characteristic
function on $C$.  By Lemma~\ref{lem:main} and induction, $F$ is
a characteristic function on $\prod_{j=0}^k C_j$ for all $n \ge k \ge 0$. 
In particular, it is a characteristic function of $C_0$ (when $k = 0$).
Theorem~\ref{main} follows.

\end{document}